\documentclass[12pt]{amsart}
\usepackage{amsmath, amsthm, amscd, amsfonts}

\newtheorem{theorem}{Theorem}[section]

\newtheorem{proposition}[theorem]{Proposition}
\newtheorem{corollary}[theorem]{Corollary}
\theoremstyle{definition}
\newtheorem{definition}[theorem]{Definition}
\newtheorem{example}[theorem]{Example}

\theoremstyle{remark}
\newtheorem{remark}[theorem]{Remark}
\numberwithin{equation}{section}
\begin{document}
\title{Approximately $C^*$-inner product preserving mappings}
\author{J. Chmieli\'nski and M. S. Moslehian}
\address{Jacek Chmieli\'nski: Institute of Mathematics, Pedagogical University of Cracow, Pod\-cho\-r\c{a}\-\.{z}ych 2, 30-084 Krak\'{o}w, Poland}
\email{jacek@ap.krakow.pl}
\address{Mohammad Sal Moslehian: Department of Mathematics, Ferdowsi University, P. O. Box 1159, Mashhad 91775, Iran}\email{moslehian@ferdowsi.um.ac.ir}
\begin{abstract}
A mapping $f: {\mathcal M} \to {\mathcal N}$ between Hilbert
$C^*$-modules approximately preserves the inner product if
\[\|\langle f(x), f(y)\rangle  - \langle x, y\rangle \| \leq \varphi(x, y),\]
for an appropriate control function $\varphi(x,y)$ and all $x, y
\in {\mathcal M}$. In this paper, we extend some results
concerning the stability of the orthogonality equation to the
framework of Hilbert $C^*$-modules on more general restricted
domains. In particular, we investigate some asymptotic behavior
and the Hyers--Ulam--Rassias stability of the orthogonality
equation.
\end{abstract}
\maketitle

\section{Introduction and preliminaries}

The notion of Hilbert $C^*$-module can be regarded as a
generalization of the concepts of Hilbert space and fibre bundle.
Hilbert $C^*$-modules were first studied by I. Kaplansky
\cite{KAP} for commutative $C^*$-algebras and later by M. A.
Rieffel \cite{RIE} and W. L. Paschke \cite{PAS} for more general
$C^*$-algebras.  These objects are useful tools in many areas such
as $AW^*$-algebra theory, theory of operator algebras, operator
K-theory, group representation theory, noncommutative geometry,
locally compact quantum groups, and theory of operator spaces; see
\cite{M-T} and references therein.

Suppose that ${\mathcal A}$ is a $C^*$-algebra and ${\mathcal M}$
is a linear space which is an algebraic left ${\mathcal
A}$-module with a compatible scalar multiplication, i.e.,
$\lambda(ax)=a(\lambda x)= (\lambda a)x$ for $x\in {\mathcal M},
a\in {\mathcal A}, \lambda\in {\mathbb C}$. The space ${\mathcal
M}$ is called a pre-Hilbert ${\mathcal A}$-module (or an inner
product ${\mathcal A}$-module) if there exists an ${\mathcal
A}$-valued inner product $\langle .,.\rangle :{\mathcal M}\times
{\mathcal M}\to {\mathcal A}$ with the following
properties :\\
(i) $\langle x, x\rangle  \geq 0$ and $\langle x, x\rangle =0$ if and only if $x=0$\\
(ii) $\langle \lambda x + y, z\rangle  = \lambda \langle x, z\rangle  + \langle y, z\rangle $\\
(iii) $\langle ax, y\rangle  = a\langle x, y\rangle $\\
(iv) $\langle x, y\rangle ^* = \langle y, x\rangle$\\
for all $x, y, z \in {\mathcal M}, a\in {\mathcal A}, \lambda\in
{\mathbb C}$. Note that the condition (i) is understood as a
statement in the $C^*$-algebra ${\mathcal A}$, where an element
$a$ is called positive if it can be represented as $bb^*$ for some
$b \in {\mathcal A}$. The conditions (ii) and (iv) implies the
inner product to be conjugate-linear in its second variable.
Validity of a useful version of the classical Cauchy-Schwartz
inequality follows that $\|x\|=\|\langle x, x\rangle
\|^\frac{1}{2}$ is a norm on ${\mathcal M}$ making it into a
normed left ${\mathcal A}$-module. The pre-Hilbert module
${\mathcal M}$ is called a Hilbert ${\mathcal A}$-module if it is
complete with respect to the above norm. Some interesting
examples are the usual Hilbert spaces as Hilbert ${\mathbb
C}$-modules, and any $C^*$-algebra ${\mathcal A}$ as a Hilbert
${\mathcal A}$-module via $\langle a,b\rangle =ab^*\;\;\;(a,b \in
{\mathcal A})$. Notice that the inner structure of a $C^*$-algebra
is essentially more complicated than complex numbers, hence the
notions such as orthogonality and theorems such as Riesz'
representation in the Hilbert space theory cannot simply be
generalized or transferred to the theory of Hilbert $C^*$-modules.

One may define an ``${\mathcal A}$-valued norm'' $|.|$ by $|x| =
\langle x, x \rangle^{1/2}$ (where, $|a|$ denotes the unique
square root of the positive element $aa^*$ in ${\mathcal A}$).
Clearly, $\|\;|x|\;\| = \|x\|$, for each $x \in {\mathcal M}$. It
is known that $|.|$ does not satisfy the triangle inequality in
general; cf. \cite{LAN}.

Roughly speaking, a functional equation $(\mathcal E)$ is {\it
stable} if any mapping which approximately satisfies the equation
$(\mathcal E)$ is near to an exact solution of $(\mathcal E)$. The
equation $(\mathcal E)$ is called {\it superstable} if any
approximate solution of $(\mathcal E)$ is, in fact, an exact
solution.

In 1940 Ulam \cite{ULA} posed the first stability problem
concerning the stability of group homomorphisms. In the next year,
Hyers \cite{HYE} gave a partial affirmative answer to the question
of Ulam in the context of Banach spaces. In 1978, Th. M. Rassias
\cite{RAS} generalized the theorem of Hyers by considering a
particular stability problem with unbounded Cauchy differences
(which is now often called the Hyers--Ulam--Rassias stability).
More general approach was considered already in 1951 by D. G.
Bourgin \cite{BOU} and later by G. L. Forti \cite{FOR}, P.
G\u{a}vruta \cite{GAV} and others. During the last decades
several stability problems for functional equations have been
investigated; we refer the reader e.g. to monographs \cite{CZE,
H-I-R1, JUN} and references therein. In particular, several
stability results have been obtained for various equations for
mappings on Hilbert $C^*$-modules, see \cite{AMY, MOS}.

A mapping $I: {\mathcal M} \to {\mathcal N}$ between Hilbert
$C^*$-modules preserves the inner product if it is a solution of
the orthogonality equation
\[\langle I(x), I(y)\rangle  = \langle x, y\rangle .\]
It is routine to show that $I$ preserves the inner product if and
only if it is ${\mathcal A}$-linear (i.e., $I(ax + \lambda y + z)
= aI(x) + \lambda I(y) + I(z)$, for all $a \in {\mathcal A}, x, y,
z \in {\mathcal M}, \lambda \in {\mathbb C}$) and it is an
isometry in the sense that $\|I(x) - I(y)\| = \|x - y\|$, for all
$x, y \in {\mathcal M}$ (for a proof in the context of Hilbert
spaces see Lemma 2.1.1 of \cite{MLA}).

A mapping $f: {\mathcal M} \to {\mathcal N}$ approximately
preserves the inner product if it satisfies
\begin{eqnarray*}
\|\langle f(x), f(y)\rangle  - \langle x, y\rangle \| \leq
\varphi(x, y),
\end{eqnarray*}
for some appropriate control function $\varphi(x,y)$ and all $x,
y \in {\mathcal M}$.

Recently, the stability of the orthogonality equation (as well as
of the so-called Wigner equation $|\langle f(x), f(y)\rangle | =
|\langle x, y\rangle |$ ) has been studied in the framework of
real and complex Hilbert spaces; see e.g. \cite{B-C, CHM, C-J}
and the Chapter (IV) of \cite{H-I-R1}. Another related work is
\cite{B-C-M} where $n$-inner product preserving mappings are
investigated.

We generalize the main results of Chmieli\'{n}ski, Badora and Jung
concerning the stability of orthogonality spaces to Hilbert
$C^*$-modules, prove the stability on a general restricted domain,
investigate some asymptotic aspects and prove the
Hyers--Ulam--Rassias stability of the orthogonality equation.

Throughout the paper, ${\mathcal M}$ and ${\mathcal N}$ denote a
pre-Hilbert module and a Hilbert module over a $C^*$-algebra
${\mathcal A}$, respectively. In addition, we denote by ${\mathbb
N}$, ${\mathbb N}_0$ and ${\mathbb R}$ the set of positive
integers, non-negative integers and real numbers, respectively. We
refer the reader to \cite{MUR} for undefined notions on
$C^*$-algebra theory and to \cite{LAN, M-T} for more information
on Hilbert $C^*$-modules.

\section{Stability on restricted domains}

Let $D$ be a subset of ${\mathcal M} \times {\mathcal M}$
containing $\Delta \times \Delta$, where $\Delta = \{x \in
{\mathcal M}: (x, x) \in D\}$, and suppose that there exists a
positive number $c \neq 1$ such that:

(i) for all $(x, y) \in D$ and all $m, n \in {\mathbb N}_0$, we
have $(c^{-n}x, c^{-m}y) \in D$;

(ii) for all $x, y \in {\mathcal M}\setminus \{0\}$ there are
nonnegative integers $m, n$ with $(c^{-n}x, c^{-m}y) \in D$.

For instance, $D$ can be chosen to be ${\mathcal M} \times
{\mathcal M}$, $\{x \in {\mathcal M}: \|x\| \leq d\} \times \{x
\in {\mathcal M}: \|x\| \leq d\}$ or $\{x \in {\mathcal M}: \|x\|
\geq d\} \times \{x \in {\mathcal M}: \|x\| \geq d\}$, where $d$
is a positive number.

Using some ideas from \cite{B-C, C-J}, we are going to extend
their main results not only to more general domains but also to a
more general framework.

\begin{theorem} \label{stab}
Consider a function $\varphi: {\mathcal M} \times {\mathcal M} \to
[0, \infty)$ satisfying
\[\lim_{m + n \to \infty} c^{m + n} \varphi(c^{-m}x, c^{-n}y) =
0,\qquad  (x, y) \in D.\] Let $f: {\mathcal M} \to {\mathcal N}$
be a mapping such that
\begin{eqnarray}\label{appinn}
\|\langle f(x), f(y)\rangle  - \langle x, y\rangle \| \leq
\varphi(x, y),\qquad (x, y) \in D.
\end{eqnarray}
Then there exist a unique ${\mathcal A}$-linear isometry $I:
{\mathcal M} \to {\mathcal N}$ and a mapping $T : {\mathcal M} \to
{\mathcal N}$ such that
\[f(x) = I(x) + T(x),\]
\[\|f(x) - I(x)\| \leq \sqrt{\varphi(x, x)},\]
\[\|T(x)\| \leq \sqrt{\varphi(x, x)},\]
\[\langle T(x), I(y)\rangle  = 0,\]
for all $x, y \in \Delta$.
\end{theorem}
\begin{proof}
For the sake of convenience, we introduce the functions $f_n:
{\mathcal M} \to {\mathcal N}$ by $f_n(x) = c^nf(c^{-n}x)$ for any
$n \in {\mathbb N}_0$. Evidently, $f_0 = f$. Recall that if $a$ is
an element of the $C^*$-algebra ${\mathcal A}$, then the real part
${\rm Re}(a)$ of $a$ is defined to be $\frac{a + a^*}{2}$. We have
also $\|{\rm Re}(a)\|\leq\|a\|$.

Let $x \in \Delta$ and $m, n \in {\mathbb N}_0$. We have
\begin{eqnarray*}
\|{\rm Re}(\langle f_n(x), f_m(x)\rangle ) - \langle x, x\rangle
\| &=& \|{\rm Re}(\langle f_n(x), f_m(x)\rangle  - \langle x,
x\rangle )\|\\
&\leq& \|\langle f_n(x), f_m(x)\rangle  - \langle x, x\rangle \|\\
&=& c^{n + m}\|\langle f(c^{-n}x), f(c^{-m}x)\rangle  - \langle c^{-n}x, c^{-m}x\rangle \|\\
&\leq& c^{n + m}\varphi(c^{-n}x, c^{-m}x),
\end{eqnarray*}
whence
\begin{eqnarray}\label{fmn}
\|f_n(x) - f_m(x)\|^2 &=& \|\;|f_n(x) - f_m(x)|\;\|^2 \nonumber\\
&=& \|\;|f_n(x) - f_m(x)|^2\;\|\nonumber\\
&=& \|\langle f_n(x) - f_m(x), f_n(x) - f_m(x)\rangle \|\nonumber\\
&\leq& \|\;|f_n(x)|^2 + |f_m(x)|^2 - 2 {\rm Re} (\langle f_n(x), f_m(x)\rangle )\|\nonumber\\
&\leq& \|\;|f_n(x)|^2 - |x|^2\| + \|\;|f_m(x)|^2 - |x|^2\|\nonumber\\
&& + 2 \|{\rm Re} (\langle f_n(x), f_m(x)\rangle ) - \langle x, x\rangle \|\nonumber\\
&\leq& c^{2n}\varphi(c^{-n}x, c^{-n}x) + c^{2m}\varphi(c^{-m}x,
c^{-m}x) + 2c^{n + m}\varphi(c^{-n}x, c^{-m}x).
\end{eqnarray}
Thus the sequence $\{f_n(x)\}$ is a Cauchy one in the complete
space ${\mathcal N}$, whence it is convergent. Set
\[I_*(x) := \lim_{n \to \infty} f_n(x),\qquad x \in \Delta.\]
Let ($x, y) \in\Delta\times\Delta\subset D$. Then
\[\|\langle f_n(x), f_n(y)\rangle  - \langle x, y\rangle \| \leq c^{2n}\varphi(c^{-n}x,
c^{-n}y),\] for all $n$. Letting $n \to \infty$ we get
\[\langle I_*(x), I_*(y)\rangle  = \langle x, y\rangle .\]
Putting $m = 0$ in (\ref{fmn}) we get
\[\|f_n(x) - f(x)\|^2 \leq c^{2n}\varphi(c^{-n}x, c^{-n}x) + \varphi(x,
x) + 2c^{n}\varphi(c^{-n}x, x)\] from which we conclude that
\begin{eqnarray} \label{I*f}
\|I_*(x) - f(x)\| \leq \sqrt{\varphi(x, x)},\qquad  x \in \Delta.
\end{eqnarray}

Let us define the mapping $I : {\mathcal M} \to {\mathcal N}$ as
\begin{eqnarray*}
I(x) := \begin{cases}
c^{n(x)}I_*(c^{-n(x)}x),& x \in {\mathcal M} \setminus \{0\};\\
0, & x = 0
\end{cases}
\end{eqnarray*}
where $n(x) = \min\{n \in {\mathbb N}_0: c^{-n}x \in \Delta\}$.
Note that if $x$ is a non-zero element in ${\mathcal M}$, then
$(c^{-n}x, c^{-m}x) \in D$ for some $n, m$. If $k = \max\{m, n\}$,
then $(c^{-k}x, c^{-k}x) \in D$ and so $c^{-k}x \in \Delta$. Hence
$I$ is well-defined. If $x \in \Delta$, then $n(x) = 0$ and so
$I(x) = I_*(x)$. It follows then from (\ref{I*f}) that
\begin{eqnarray}\label{If}
\|I(x) - f(x)\| \leq \sqrt{\varphi(x, x)},\qquad x \in \Delta.
\end{eqnarray}

We are going to prove that $I$ is an inner product preserving
mapping and so it is an isometry. To see this, assume that $x, y
\in {\mathcal M}$. If $x = 0$ or $y = 0$, then $\langle I(x),
I(y)\rangle  = 0 = \langle x, y\rangle $. Let $x \neq 0$ and $y
\neq 0$. Then
\begin{eqnarray*}
\langle I(x), I(y)\rangle  &=& \langle c^{n(x)}I_*(c^{-n(x)}x), c^{n(y)}I_*(c^{-n(y)}y)\rangle \\
&=& c^{n(x) + n(y)} \langle I_*(c^{-n(x)}x), I_*(c^{-n(y)}y)\rangle \\
&=& c^{n(x) + n(y)} \langle c^{-n(x)}x, c^{-n(y)}y\rangle \\
&=& \langle x, y\rangle .
\end{eqnarray*}
For proving the uniqueness assertion, consider inner product
preserving mappings $I_1, I_2$ satisfying $\|I_j(x) - f(x)\| \leq
\sqrt{\varphi(x, x)}\;\; (j = 0, 1)$ for all $x \in \Delta$. First
note that for each $x \in \Delta$ and all $n\in{\mathbb N}_0$ we
have
\begin{eqnarray*}
\|I_1(x) - I_2(x)\| &=& c^n\|I_1(c^{-n}x) - I_2(c^{-n}x)\|\\
&\leq& c^n\|I_1(c^{-n}x) - f(c^{-n}x)\| + c^n\|I_2(c^{-n}x) -
f(c^{-n}x)\| \\
&\leq& 2 \sqrt{c^{2n}\varphi(c^{-n}x, c^{-n}x)},
\end{eqnarray*}
whence $I_1(x) = I_2(x)$ on $\Delta$. Now for each $x \in
{\mathcal M}$, there exists $n(x) \in {\mathbb N}_0$ such that
$c^{-n(x)}x \in \Delta$. Therefore
\[I_1(x) = c^{n(x)}I_1(c^{-n(x)}x) = c^{n(x)}I_2(c^{-n(x)}x) =
I_2(x).\] Next, put $T(x) = f(x) - I(x)$. Then (\ref{If}) yields
$\|T(x)\| \leq \sqrt{\varphi(x, x)}$ for all $x \in \Delta$.

Let $(x, y) \in D$, then $(x,c^{-n(y)}y)\in D$ and $c^{-n(y)}y\in
\Delta$. Then $(x, c^{-n}c^{-n(y)}y) \in D$ for all $n$. Therefore
(\ref{appinn}) yields
\[\|\langle f(x), f_n(c^{-n(y)}y)\rangle  - \langle x, c^{-n(y)}y\rangle \| \leq c^n\varphi(x, c^{-n}c^{-n(y)}y).\]
Thus
\[\langle f(x), I_*(c^{-n(y)}y)\rangle = \langle
x,c^{-n(y)}y\rangle,\] whence $\langle f(x), I(y)\rangle  =
\langle x, y\rangle $, and
\[\langle T(x), I(y)\rangle  = \langle f(x) - I(x), I(y)\rangle = \langle f(x), I(y)\rangle
- \langle I(x), I(y)\rangle  = \langle f(x), I(y)\rangle  -
\langle x, y\rangle = 0.\]
\end{proof}
\begin{remark}
If $f$ is a function such that $f(cx) = cf(x)$, then $f(0) = 0$
and $I_*(x) = \lim_{n \to \infty} c^nf(c^{-n}x) = f(x)$ for all $x
\in \Delta$. It follows that $f(x) = I(x)$ for all $x \in
{\mathcal M}$.
\end{remark}
The following example, which is a slight modification of Example 1
of \cite{CHM}, shows that the bound $\sqrt{\varphi(x, y)}$ in
(\ref{If}) is sharp and we have no control on the bounded function
$T$. This means that $T$ is neither additive nor continuous in
general.
\begin{example}
Let ${\mathcal M}, {\mathcal N}$ be the Hilbert space $\ell^2$.
Assume that $g: {\mathcal M} \to {\mathbb C}$ is an arbitrary
mapping satisfying $|g(x)| \leq \sqrt{\varphi(x, y)}$. Define the
mapping $f: {\mathcal M} \to {\mathcal N}$ by $f(x) = (g(x), t_1,
t_2, \ldots)$ where $x = (t_1, t_2, \ldots) \in {\mathcal M}$.
Clearly, $\langle f(x), f(y)\rangle  = |g(x)|^2 + \langle x,
y\rangle $, for all $x, y \in {\mathcal M}$. Then $I((t_1, t_2,
\ldots)) = (0, t_1, t_2, \ldots)$ and $T(x) = (g(x), 0, 0,
\ldots)$ are the unique mappings fulfilling the required
conditions in Theorem \ref{stab}.
\end{example}

\begin{corollary}
Suppose that either $p, q > 1$ or $p, q < 1$  are real numbers and
$\alpha > 0$. Let $f: {\mathcal M} \to {\mathcal N}$ be a mapping
such that
\begin{eqnarray*}
\|\langle f(x), f(y)\rangle  - \langle x, y\rangle \| \leq \alpha
\|x\|^p\|y\|^q,\qquad \mbox{for all}\ x,y\in{\mathcal M}.
\end{eqnarray*}
Then there exists a unique linear isometry $I: {\mathcal M} \to
{\mathcal N}$ such that
\[\|f(x) - I(x)\| \leq \sqrt{\alpha}\|x\|^{\frac{p + q}{2}},\]
for all $x \in {\mathcal M}$.
\end{corollary}
\begin{proof}
Let $\varphi(x, y) = \alpha \|x\|^p\|y\|^q$. Consider $D =
{\mathcal M} \times {\mathcal M}$ together with $c >1$ if $p, q
>1$; and $c < 1$ if $p, q < 1$.
\end{proof}
\begin{remark}
The above result holds true also in cases $p=1,\ q\neq 1$ or
$p\neq 1,\ q=1$. The Corollary is not true for $p =q= 1$, in
general. For a counterexample see Example 2 of \cite{B-C}.
\end{remark}

In a particular case, where ${\mathcal M}$ and ${\mathcal N}$ are
of the same finite dimension we can prove superstability.

\begin{proposition}
Let $\dim{\mathcal M}=\dim{\mathcal N}<\infty$. Suppose that
$f:{\mathcal M}\to {\mathcal N}$ satisfies (\ref{appinn}) with
$\varphi$ as in Theorem \ref{stab}. Then there exists a linear
isometry $I:{\mathcal M}\to{\mathcal N}$ such that $f=I$ on
$\Delta$.
\end{proposition}
\begin{proof}
Let $I:{\mathcal M}\to{\mathcal N}$ be the linear isometry from
the assertion of Theorem \ref{stab} and $T=f-I$. $I$ maps
${\mathcal M}$ onto a subspace $I({\mathcal M})$ of ${\mathcal
N}$. Since $\dim{\mathcal M}=\dim I({\mathcal M})$, and ${\mathcal
M}$ and ${\mathcal N}$ are of the same finite dimension, we get
$I({\mathcal M})={\mathcal N}$. For $x\in \Delta$ we have
$T(x)\bot I({\mathcal M})$, i.e., $T(x)\bot {\mathcal N}$ whence
$T(x)=0$. Thus $f=I$ on $\Delta$.
\end{proof}
Taking $D={\mathcal M}\times{\mathcal M}$ we get immediately:
\begin{corollary}
Let $\dim{\mathcal M}=\dim{\mathcal N}<\infty$ and suppose that
$f:{\mathcal M}\to {\mathcal N}$ satisfies
\begin{eqnarray}\label{appinn2}
\|\langle f(x), f(y)\rangle  - \langle x, y\rangle \| \leq
\varphi(x, y),\qquad x, y \in {\mathcal M}
\end{eqnarray}
where $\varphi: {\mathcal M} \times {\mathcal M} \to [0, \infty)$
satisfies (with some $0<c\neq 1$)
\[\lim_{m + n \to \infty} c^{m + n} \varphi(c^{-m}x, c^{-n}y) =
0,\qquad \mbox{for all}\  x, y \in {\mathcal M}.
\]
Then $f$ is an inner product preserving mapping.
\end{corollary}

\section{Asymptotic behavior of orthogonality equation}

Following \cite{H-I-R2}, a mapping $f: {\mathcal M} \to {\mathcal
N}$ is called \emph{ $p$-asymptotically close to an isometry
mapping} $I$ if $\lim_{\|x\| \to \infty}\frac{\|f(x) -
I(x)\|}{\|x\|^p} = 0$.
\begin{definition}
A mapping $f: {\mathcal M} \to {\mathcal N}$ satisfies
\emph{$p$-asymptotically} the orthogonality equation if for each
$\varepsilon > 0$ there exists $K > 0$ such that
\begin{eqnarray}\label{assipp}
\|\langle f(x), f(y)\rangle - \langle x, y\rangle \| \leq
\varepsilon \|x\|^p\;\|y\|^p,
\end{eqnarray}
for all $x, y \in {\mathcal M}$ such that $\max \{\|x\|, \|y\|\}
\geq K$.
\end{definition}

We are going to show that asymptotically orthogonality preserving
mappings are asymptotically close to isometries.

\begin{theorem}
If $0 < p < 1$ and a mapping $f: {\mathcal M} \to {\mathcal N}$
satisfies $p$-asymptotically the orthogonality equation, then it
is $p$-asymptotically close to a linear isometry mapping.
\end{theorem}
\begin{proof}
By the assumption $f$ satisfies $p$-asymptotically the
orthogonality equation, hence there exists $K_0 > 0$ such that
\[\|\langle f(x), f(y)\rangle - \langle x, y\rangle \| \leq
\|x\|^p\;\|y\|^p\] for all $x, y \in {\mathcal M}$ with $\max
\{\|x\|, \|y\|\} \geq K_0$. It follows from Theorem \ref{stab}
(for $D=\{x:\ \|x\|\geq K_0\}\times{\mathcal M}\cup {\mathcal
M}\times\{x:\ \|x\|\geq K_0\}$, $\Delta=\{x:\ \|x\|\geq K_0\}$,
$0<c<1$ and $\varphi(x,y):=\|x\|^p\|y\|^p$) that there exists a
linear isometry $I_0$ such that
\begin{eqnarray}\label{fI0}
\|f(x) - I_0(x)\| \leq \|x\|^p
\end{eqnarray}
for all $x$ with $\|x\| \geq K_0$.

Given $\varepsilon > 0$, the assumption gives again a number
$K_\varepsilon \geq K_0$ such that
\[\|\langle f(x), f(y)\rangle - \langle x, y\rangle \| \leq
\varepsilon \|x\|^p\;\|y\|^p,\] for all $x, y \in {\mathcal M}$
with $\max \{\|x\|, \|y\|\} \geq K_\varepsilon$. Applying again
Theorem \ref{stab} we get an isometry $I_\varepsilon$ such that
\begin{eqnarray}\label{fIepsilon}
\|f(x) - I_\varepsilon(x)\| \leq \sqrt{\varepsilon}\|x\|^p
\end{eqnarray}
for all $x$ with $\|x\| \geq K_\varepsilon$.

We claim that $I_\varepsilon = I_0$. To see this, let $x \in
{\mathcal M}\setminus\{0\}$ be an arbitrary element. There exists
$N$ such that for all $n > N$, $\|2^nx\| \geq K_\varepsilon \geq
K_0$. By (\ref{fI0}) and (\ref{fIepsilon}) we have
\begin{eqnarray*}
\|I_\varepsilon(x) - I_0(x)\| &=& 2^{-n}\|I_\varepsilon(2^nx) -
I_0(2^nx)\|\\
 &\leq& 2^{-n}\|I_\varepsilon(2^nx) - f(2^nx)\| +
2^{-n}\|f(2^nx) - I_0(2^nx)\|\\
&\leq& 2^{(p - 1)n}(\sqrt{\varepsilon} + 1)\|x\|^p.
\end{eqnarray*}
The right hand side tends to zero as $n \to \infty$, hence
$I_\varepsilon = I_0$. Thus (\ref{fIepsilon}) implies that
\[\frac{\|f(x) - I_0(x)\|}{\|x\|^p} < \sqrt{\varepsilon}\]
for all $x$ with $\|x\| \geq K_\varepsilon$. Thus $f$ is
$p$-asymptotically close to the isometry mapping $I_0$.\end{proof}

\begin{remark}
Assume that $p>1$ and $f:{\mathcal M}\to {\mathcal N}$ is such
that for each $\varepsilon>0$ there exists $K>0$ such that
(\ref{assipp}) holds for all $x,y\in{\mathcal M}$ satisfying
$\min\{\|x\|,\|y\|\}\leq K$. Analogously as above, one can prove
that there exists a linear  isometry $I:{\mathcal M}\to {\mathcal
N}$ such that
\begin{eqnarray*}
\lim_{\|x\|\to 0}\frac{\|f(x)-I(x)\|}{\|x\|^p}=0.
\end{eqnarray*}
\end{remark}

\section{Hyers--Ulam--Rassias stability}

In this section, we prove the Hyers--Ulam--Rassias stability of
the orthogonality equation.

\begin{theorem} \label{HUR-stab}
Let $f: {\mathcal M} \to {\mathcal N}$ be an approximately inner
product preserving mapping on ${\mathcal M}$ associated with a
control function $\varphi: {\mathcal M} \times {\mathcal M} \to
[0, \infty)$. We assume that the control function $\psi$ defined
by
\begin{eqnarray*}
\psi(x, y) &:=& \big(\varphi(x + y, x + y) + \varphi(x, x + y) +
\varphi(y, x + y) + \varphi(x + y, x)\\
&& + \varphi(x, x) + \varphi(y, x) + \varphi(x + y, y) +
\varphi(x, y) + \varphi(y, y)\big )^{1/2}
\end{eqnarray*}
satisfies either
\begin{eqnarray}\label{1/2}
\widetilde{\psi}(x): =\sum_{n=0}^{\infty} 2^{-n-1} \psi(2^nx,
2^nx) < \infty,
\end{eqnarray}
or
\begin{eqnarray}\label{2}
\widetilde{\psi}(x): =\sum_{n=1}^{\infty} 2^{n - 1} \psi(2^{-n}x,
2^{-n}x) < \infty
\end{eqnarray}
for all $x\in {\mathcal M}$. Then there exists a unique linear
isometry $I: {\mathcal M} \to {\mathcal N}$ such that
\[\|f(x) - I(x)\| \leq \widetilde{\psi}(x, x),\]
\end{theorem}
\begin{proof}
Let $x, y, z \in {\mathcal M}$ and put $A = f(x + y) - f(x) -
f(y)$ . We have
\begin{eqnarray*}
\|\langle A, f(z)\rangle \| &\leq& \|\langle f(x + y),
f(z)\rangle  - \langle x
+ y, z\rangle \|\\
&& + \|\langle f(x), f(z)\rangle  -\langle x, z\rangle \| + \|\langle f(y), f(z)\rangle  - \langle y, z\rangle \|\\
&\leq& \varphi(x + y, z) + \varphi(x, z) + \varphi(y, z),
\end{eqnarray*}
whence
\begin{eqnarray*}
\|f(x +y) - f(x) - f(y)\|^2 &=& \|\langle A, f(x + y) - f(x) - f(y)\rangle \|\\
&\leq& \|\langle A, f(x + y)\rangle \| + \|\langle A, f(x)\rangle \| +\|\langle A, f(y)\rangle \|\\
&\leq& \varphi(x + y, x + y) + \varphi(x, x + y) + \varphi(y, x +
y)+ \varphi(x + y, x)\\
&& + \varphi(x, x) + \varphi(y, x) + \varphi(x + y, y) +
\varphi(x, y) + \varphi(y, y).
\end{eqnarray*}
It follows that
\begin{eqnarray*}
\|f(x +y) - f(x) - f(y)\| \leq  \psi(x, y),
\end{eqnarray*}
whence, in particular,
\[
\|f(2x)-2f(x)\|\leq\psi(x,x),\qquad x\in{\mathcal M}.
\]
Using the induction, one can easily verify the following
inequalities:
\begin{eqnarray}\label{Cau1}
\|2^{-n}f(2^n x) - 2^{-m}f(2^m x)\| \leq
\sum_{k=m}^{n-1}2^{-k-1}\psi(2^kx, 2^kx),
\end{eqnarray}
\begin{eqnarray}\label{Cau2}
\|2^nf(2^{-n} x) - 2^mf(2^{-m} x)\| \leq \sum_{k= m + 1
}^{n}2^{k-1}\psi(2^{-k}x, 2^{-k}x)
\end{eqnarray}
for all integers $n > m \geq 0$ and $x \in {\mathcal M}$. It
follows that the sequence $\{c^nf(c^{-n}x)\}$ with $c=\frac{1}{2}$
or $c=2$, respectively, is a Cauchy one, whence it is convergent.
Define the mapping $I: {\mathcal M} \to {\mathcal N}$ by $I(x) :=
\lim_{n \to \infty} c^nf(c^{-n})$. Since $f$ is approximately
inner product preserving, we have
\begin{eqnarray*}
c^{2n}\|\langle f(c^{-n}x), f(c^{-n}y)\rangle  - \langle c^{-n}x,
c^{-n}y\rangle \| \leq c^{2n}\varphi(c^{-n}x, c^{-n}y).
\end{eqnarray*}
Passing to the limit as $n$ tends to infinity we get
\[
\langle I(x), I(y)\rangle  = \langle x, y\rangle,\qquad x, y \in
{\mathcal M}.
\]
In addition, it follows from (\ref{Cau1}) and (\ref{Cau2}) with
$m = 0$ as $n \to \infty$ that
\[\|f(x) - I(x)\| \leq \widetilde{\psi}(x, x),\]
\end{proof}

\begin{corollary}
Suppose that $p \neq 2$ is a real number and $\beta > 0$. Let $f:
{\mathcal M} \to {\mathcal N}$ be a mapping such that
\begin{eqnarray*}
\|\langle f(x), f(y)\rangle  - \langle x, y\rangle \| \leq \beta
(\|x\|^p + \|y\|^p),
\end{eqnarray*}
for all $x, y \in {\mathcal M}$. Then there exists a unique linear
isometry $I: {\mathcal M} \to {\mathcal N}$ such that
\[\|f(x) - I(x)\| \leq \frac{\sqrt{6\beta (2^p + 2)}}{\left|2^{\frac{p}{2}} -
2\right|}\,\|x\|^{\frac{p}{2}},\qquad  \mbox{for all}\ x \in
{\mathcal M}.
\]
\end{corollary}
\begin{proof}
Apply Theorem \ref{HUR-stab} with $\varphi(x, y) = \beta (\|x\|^p
+ \|y\|^p)$ and consider (\ref{1/2}) if $p < 2$, and (\ref{2}) if
$p > 2$.
\end{proof}
\begin{remark}
The case $p = 2$ remains unsolved.
\end{remark}


\begin{thebibliography}{99}
\bibitem{AMY} M. Amyari, \textit{Stability of $C^*$-inner products}, to appear in J. Math. Anal. Appl.
\bibitem{B-C-M} C. Baak, H. Chu and M. S. Moslehian, \textit{On the Cauchy--Rassias inequality and linear $n$-inner product preserving
mappings}, to appear in Math. Inequ. Appl., arXiv:
math.FA/0501159.
\bibitem{B-C} R. Badora and J. Chmieli\'{n}ski, \textit{Decomposition of mappings approximately inner
product preserving}, Nonlinear Anal. (TMA) 62 (2005) 1015-–1023.
\bibitem{BOU} D. G. Bourgin, \textit{Classes of transformations and bordering
transformations}, Bull. Amer. Math. Soc.  57 (1951), 223--237.
\bibitem{CHM} J. Chmieli\'{n}ski, {\it On a singular case in the Hyers--Ulam--Rassias stability of the Wigner equation}, J. Math.
Anal. Appl. 289 (2004), 571--583.
\bibitem{C-J} J. Chmieli\'{n}ski and S.-M. Jung, \textit{On the stability of the Wigner equation on a restricted domain}, J. Math. Anal.
Appl. 254 (2001) 309–-320.
\bibitem{CZE} S. Czerwik, \textit{Functional Equations and Inequalities in Several Variables}, World Scientific, River Edge, NJ, 2002.
\bibitem{FOR} G. L. Forti, \textit{An existence and stability theorem for a class
of functional equations}, Stochastica 4 (1980), 23--30.
\bibitem{GAV} P. G\u{a}vruta, \textit{A generalization of the Hyers--Ulam--Rassias stability of approximately additive mappings},
J. Math. Anal. Appl., 184 (1994), 431--436.
\bibitem{HYE} D. H. Hyers, \textit{On the stability of the linear functional equation}, Proc. Nat. Acad. Sci. U.S.A. 27 (1941), 222--224.
\bibitem{H-I-R1} D. H. Hyers, G. Isac and Th. M. Rassias, \textit{Stability of Functional Equations in Several Variables}, Birkh\" auser, Basel, 1998.
\bibitem{H-I-R2} D. H. Hyers, G. Isac and Th. M. Rassias, \textit{On the asymptoticity aspect of
Hyers–Ulam stability of mappings}, Proc. Amer. Math. Soc. 126
(1998), 425–-430.
\bibitem{JUN} S.-M. Jung, \textit{Hyers--Ulam--Rassias Stability of Functional Equations in Mathematical Analysis},
Hadronic Press, 2001.
\bibitem{KAP} I. Kaplansky, \textit{Modules over operator algebras}, Amer. J. Math. 75 (1953), 839--858.
\bibitem{LAN} E. C. Lance, \textit{Hilbert $C^*$-Modules}, LMS Lecture Note Series 210, Cambridge Univ. Press, 1995.
\bibitem{M-T} V. M. Manuilov and E. V. Troitsky, \textit{Hilbert
$C^*$-modules}, Translations of Mathematical Monographs, 226.
American Mathematical Society, Providence, RI, 2005
\bibitem{MLA} W. Mlak, \textit{Hilbert Spaces and Operator Theory}, Kluwer Academic Publishers, PWN-Polish Scientific Publishers,
Dodrecht,Warszawa, 1991.
\bibitem{MOS} M. S. Moslehian, \textit{Stability of adjointable mappings in Hilbert
$C^*$-modules}, arXiv: math.FA/0501139.
\bibitem{MUR} J. G. Murphy, \textit{Operator Theory and $C^*$-Algebras}, Acad. Press, 1990.
\bibitem{PAS} W. L. Paschke, \textit{Inner product modules over $B^*$-algebras}, Trans. Amer. Math. Soc. 182 (1973),443--468.
\bibitem{RAS} Th. M. Rassias, \textit{On the stability of the linear mapping in Banach spaces}, Proc. Amer. Math. Soc. 72 (1978), 297--300.
\bibitem{RIE} M. A. Rieffel, \textit{Morita equivalence representations of $C^*$-algebras}, Adv. in Math. 13 (1974), 176--257.
\bibitem{ULA} S. M. Ulam, \textit{Problems in Modern Mathematics}, Chapter VI, Science Editions, Wiley, New York, 1964.
\end{thebibliography}
\end{document}